\documentclass[a4paper]{amsart}
\usepackage{epsfig,psfrag,isolatin1,amsmath,amssymb}
\newtheorem{theorem}{Theorem}
\newtheorem{proposition}{Proposition}
\newtheorem{lemma}{Lemma}

\theoremstyle{definition}
\newtheorem{definition}{Definition}

\theoremstyle{remark}

\newcommand{\jj}[6]{\left|\begin{smallmatrix} #1 & #2 & #3\\
                                              #6 & #5 & #4
                    \end{smallmatrix}\right|}
\def\co{\colon\thinspace}

\newcommand{\eps}{\epsilon}

\renewcommand{\d}{\partial}

\newcommand{\begriff}[1]{\textbf{#1}}                  

\begin{document}
\author{Simon A. King}\thanks{Partially supported by INTAS project CALCOMET-GT ref.~03-51-3663.}
\title{Ideal Turaev--Viro invariants}  
\date{September 8, 2005}

\maketitle
\begin{abstract}
  A Turaev--Viro invariant is a state sum, i.e., a polynomial that can
  be read off from a special spine or a triangulation of a compact
  $3$-manifold. If the polynomial is evaluated at the solution of a
  certain system of polynomial equations (Biedenharn--Elliott
  equations) then the result is a homeomorphism invariant of the
  manifold (``numerical Turaev-Viro invariant''). The equation system
  defines an ideal, and actually the coset of the polynomial with
  respect to that ideal is a homeomorphism invariant as well (``ideal
  Turaev--Viro invariant'').
  
  It is clear that ideal Turaev--Viro invariants are at least as strong
  as numerical Turaev--Viro invariants, and we show that there is
  reason to expect that they are strictly stronger. They offer a more
  unified approach, since many numerical Turaev--Viro invariants can be
  captured in a singly ideal Turaev--Viro invariant. Using computer
  algebra, we obtain computational results on some examples of
  ideal Turaev--Viro invariants.
\end{abstract}

\section{Introduction}
\label{sec:Intro}

Let $M_1$, $M_2$ be compact $3$-manifolds, represented by special
spines $P_1$, $P_2$ (Definition~\ref{def:Spines}). If $M_1$ is
homeomorphic to $M_2$, then $P_1$ and $P_2$ are related by a finite
sequence of certain local transformations. Turaev--Viro invariants,
originally formulated for triangulations rather than special
spines~\cite{TuraevViro92}, can be read off from any special spine of
a compact $3$-manifold: One computes the state sum, i.e. a polynomial
whose summands correspond to different ``colourings'' of the special
spine. This polynomial is a homeomorphism invariant, provided its
variables satisfy the so-called Biedenharn--Elliott equations known
from quantum physics~\cite{landau}. It is difficult to find solutions
of the Biedenharn--Elliott equations, but an important class of
solutions is provided by the representation theory of Quantum
Groups~\cite{turaev}. We call this a ``numerical Turaev--Viro
invariant''.

The starting point of this paper is the observation that the coset of
the state sum with respect to the ideal generated by the
Biedenharn--Elliott equations is a homeomorphism invariant of compact
$3$-manifolds. Hence it is not needed to find explicit solutions of
the equations. We call this an ``ideal Turaev--Viro invariant''. We
also define an invariant that captures \emph{all} numerical invariants
obtained from evaluating an ideal Turaev--Viro invariant, without the
need to compute any explicit solution. We give reasons why one should
expect that ideal Turaev--Viro invariants are strictly stronger than
all their associated numerical invariants together. Using software for
the computation of Gröbner bases, we computed some examples of ideal
Turaev--Viro invariants for closed orientable irreducible manifolds of
complexity up to $9$.

The paper is organised as follows. In Section~\ref{sec:Spines}, we
recall some basic facts about special spines of compact $3$-manifolds.
In Section~\ref{sec:statesums} we define the state sum associated to a
special $2$-polyhedron. In Section~\ref{sec:TVinvariants} we define
ideal and numerical Turaev--Viro invariants, construct an invariant
that captures all numerical Turaev--Viro invariant associated to an
ideal invariant, and observe a lower bound for the complexity of a
manifold in terms of the ideal Turaev-Viro invariant (unfortunately,
the bound turned out to be trivial in all examples that we computed).
In Section~\ref{sec:implementation} we come to the problem of how to
explicitly compute ideal Turaev--Viro invariants, based on
implemented algorithms of commutative algebra.
The computation of ideal Turaev--Viro invariants can be pretty
complex, hence it seems reasonable to introduce simplifying
assumptions. We suggest different types of simplification in
Section~\ref{sec:assumptions}.
In Section~\ref{sec:expl}, we present four examples of ideal
Turaev--Viro invariants in more detail.
In the final Section~\ref{sec:computations} we report our
computational results on these examples.

\section{Special spines of compact $3$-manifolds}
\label{sec:Spines}

\begin{definition}\label{def:Spines}
  A \begriff{simple 2-polyhedron} $P$ is a compact connected hausdorff
  space such that any point has an open neighbourhood of one of the
  following three homeomorphism types (where the point under
  consideration is marked by a thick dot):
  \\[2mm]

  \centerline{{\psfrag{i}{\emph{(i)}}\psfrag{ii}{\emph{(ii)}}\psfrag{iii}{\emph{(iii)}}
      \epsfig{width=0.9\linewidth,file=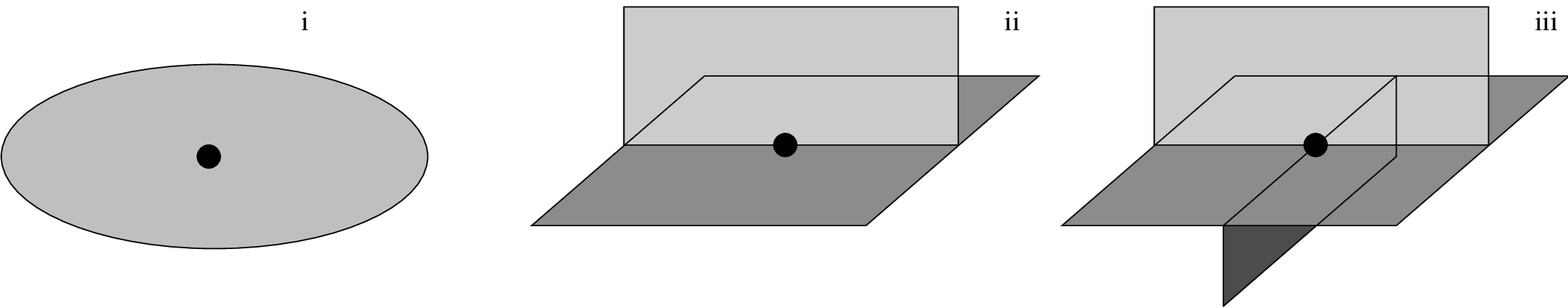}}}
  The connected components of points of type \emph{(i)} are the
  \begriff{2-strata} of $P$, the connected components of points of
  type \emph{(ii)} are the \begriff{true edges} of $P$, and the points
  of type \emph{(iii)} are the \begriff{true vertices} of $P$. The
  4-valent graph $S(P)\subset P$ formed by true edges and true
  vertices is the \begriff{singular graph} of $P$. The set of 2-strata
  of $P$ is denoted by $\mathcal C(P)$, the set of true edges of $P$
  is denoted by $\mathcal E(P)$, and the set of true vertices of $P$
  is denoted by $\mathcal V(P)$. 
  
  A simple 2-polyhedron is special, if it has a true vertex, its
  singular graph is connected, and its 2-strata are homeomorphic to
  open discs. Let $M$ be a compact 3-manifold. A special 2-polyhedron
  $P$ embedded in $M$ is a \begriff{special spine} of $M$, if $\d
  M=\emptyset$ and $M\setminus P$ is homeomorphic to a 3-ball, or if
  $\d M\not=\emptyset$ and $M\setminus P\approx (\d M)\times [0,1)$
  (where $[0,1)$ denotes a half-open interval).\end{definition}

A general reference for the theory of special spines of compact
$3$-manifolds is~\cite{matveev}. Any compact 3-manifold has a special
spine, which can be deduced from the fact that any compact 3-manifold
admits a triangulation~\cite{moise1}. Moreover, the homeomorphism type
of a special spine uniquely determines the homeomorphism type of the
3-manifold~\cite{Casler65}. The following classical result explains
how all special spines of a compact 3-manifold are related with each
other.

\begin{figure}[htbp]
  \centering
  {\psfrag{T}{$T$}\psfrag{T-}{$T^{-1}$}\epsfig{width=\linewidth,file=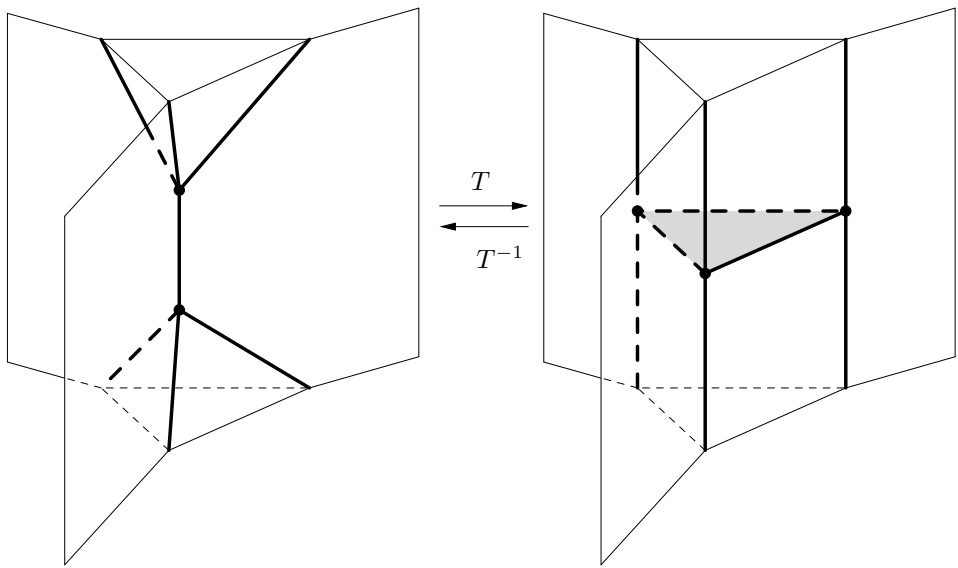}}  
  \caption{The Matveev-Piergallini move}
  \label{fig:move}
\end{figure}

\begin{theorem}[Matveev~\cite{Matveev88}, Piergallini~\cite{Piergallini88}]
  \label{thm:matvpierg} 
  Let $M$ be a compact 3-manifold with special spines $P_1,P_2$, and
  assume that $P_1$ and $P_2$ both have at least two true vertices.
  Then $P_1$ and $P_2$ are related by a finite sequence of a local
  transformation called $T$ move and its inverse. The $T$ move is
  shown in Figure~\ref{fig:move}, where true vertices are marked by a
  thick dot and true edges are drawn bold.
\end{theorem}

``Local tranformation'' means that the special 2-polyhedron remains
unchanged outside of the depicted part. It is remarkable that a single
type of local transformation suffices. When working with
triangulations of 3-manifolds, there is a similar transformation
result at hand, due to Pachner~\cite{pachner} --- but this one uses
two different types of local transformations.
Theorem~\ref{thm:matvpierg} offers a strategy for constructing
homeomorphism invariants of 3-manifolds: Define some algebraic
expression that can be read off from any special 2-polyhedron with at
least two true vertices, study how this expression changes under the
moves $T^\pm$, and quotient out these changes.

\section{Turaev--Viro state sums}
\label{sec:statesums}

Let $P$ be a special 2-polyhedron. Let $\mathcal F$ be a finite set,
to whose elements we will refer by \begriff{2-strata colours}. A
\begriff{$\mathcal F$-colouring} of $P$ is any map $\varphi\co
\mathcal C(P)\to \mathcal F$, i.e., to any 2-stratum of $P$ is
assigned a colour. Hence we can consider $\varphi$ as a locally
constant map $P\setminus S(P)\to \mathcal F$, and in this sense it is
clear what we mean by the restriction of $\varphi$ to a subset of $P$.
By $\Phi_{\mathcal F}(P)$ we denote the set of all $\mathcal
F$-colourings of $P$.

Let $\mathcal G$ be another finite set, to whose elements we will
refer by \begriff{edge colours}. A \begriff{$\mathcal{F,G}$-colouring}
of $P$ is any pair $(\varphi,\psi)$ of maps $\varphi\co \mathcal
C(P)\to \mathcal F$, $\psi\co \mathcal E(P)\to \mathcal G$. Again, we
can consider $\varphi$ as a locally constant map $P\setminus S(P)\to
\mathcal F$ and $\psi$ as a locally constant map $S(P)\setminus
\mathcal V(P)\to \mathcal G$, and in this sense we can talk about the
restriction of $(\varphi,\psi)$ to a subset of $P$.  We denote by
$\Phi_{\mathcal{F,G}}(P)$ the set of all $\mathcal{F,G}$-colourings of
$P$. 

\begin{figure}[htbp]
  \centering
  {\psfrag{a}{$a$}\psfrag{b}{$b$}\psfrag{c}{$c$}\psfrag{d}{$d$}\psfrag{e}{$e$}\psfrag{f}{$f$}
 \psfrag{A}{$A$}\psfrag{B}{$B$}\psfrag{C}{$C$}\psfrag{D}{$D$}\epsfig{scale=0.8,file=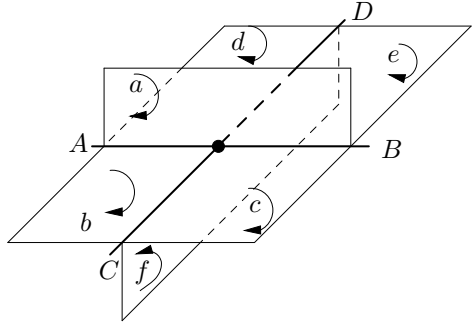}}
  \caption{A true vertex $v$ with coloured neighbourhood}
  \label{fig:6j4k}
\end{figure}

\begin{figure}[htbp]
  \centering
  {\psfrag{a}{$a$}\psfrag{b}{$b$}\psfrag{c}{$c$}\psfrag{d}{$d$}\psfrag{e}{$e$}\psfrag{f}{$f$}
 \psfrag{A}{$A$}\psfrag{B}{$B$}\psfrag{C}{$C$}\psfrag{D}{$D$}
  \epsfig{width=\linewidth,file=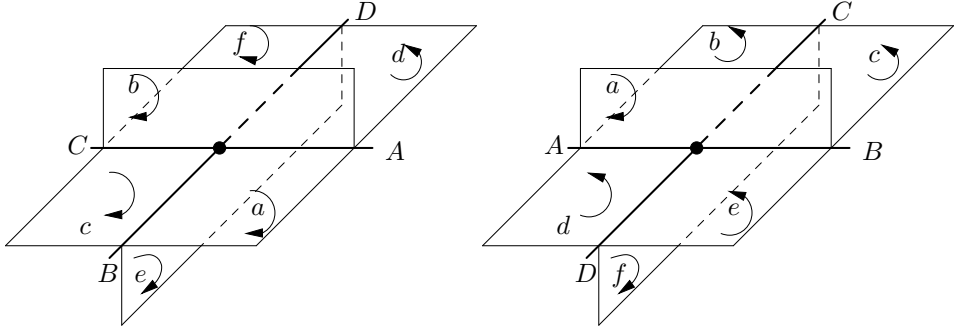}}
  \caption{Symmetry of $6j4k$-symbols}
  \label{fig:symmetry}
\end{figure}

Let us now choose an orientation independently for each 2-stratum of
$P$. We would like to make the colour of a 2-stratum dependent on the
orientation. For that purpose, we take an involution ``$-$'' on
$\mathcal F$, and impose that, if an oriented 2-stratum is
couloured by $f\in\mathcal F$, then the oppositely oriented 2-stratum has
the colour $-f\in\mathcal F$. For any $f\in \mathcal F$, the
symbol $w(f)$ is referred to as the \begriff{weight} of $f$.  At a
true vertex of $P$, six 2-strata and four true edges meet (counted
with multiplicities). A $\mathcal{F,G}$-colouring $(\varphi,\psi)$ of
$P$ thus yields for each true vertex of $P$ a 6-tuple of 2-strata
coulours together with a 4-tuple of edge coulours.  If to the 2-strata
and true edges in the neighbourhood of a true vertex $v$ are assigned
by $\varphi$ and $\psi$ coulours $a,\ldots,f\in\mathcal F$ and
$A,\ldots,D\in\mathcal G$, respectively, as depicted in
Figure~\ref{fig:6j4k} (where circular orientations of the 2-strata are
indicated by arrows, the true vertex is marked by a thick dot, and the
true edges are drawn bold),
we associate to $v$ a symbol $v^{\varphi,\psi} := \jj
abcdef^{A,B}_{C,D}$, called \begriff{$6j4k$-symbol}. If essentially
only the 2-strata are coloured, i.e., if $\mathcal G=\{*\}$ contains a
single colour, we just have a \begriff{$6j$-symbol} $v^\varphi:=\jj
abcdef=\jj abcdef^{*,*}_{*,*}$, to simplify the notation.  We shall
need that the weights and $6j4k$-symbols only depend on the
colourings, but not on additional choices. Therefore the weight of the
colour of a 2-stratum must not depend on the choice of orientation of
the 2-stratum, hence we assume that $w(f)=w(-f)$ for all $f\in\mathcal
F$.  When depicting a true vertex $v$ as in Figure~\ref{fig:6j4k},
this also involves some choices: Which of the six 2-strata is on top?
And which side of the top 2-stratum shall be visible in front? By the
tetrahedral symmetry of true vertices, the same vertex could also be
depicted as in Figure~\ref{fig:symmetry}.
Our $6j4k$-symbols (resp. $6j$-symbols) shall have the same symmetry.
Hence, we assume that the following identities hold for all
$a,b,c,d,e,f \in \mathcal F$ and all $A,B,C,D\in \mathcal G$:
\begin{eqnarray*}
  \jj abcdef^{A,B}_{C,D}  &=& \jj bcaf{-d}{-e}^{C,A}_{B,D}
                \\ 
        &=& \jj a{-d}{-e}{-b}{-c}{-f}^{A,B}_{D,C}
\end{eqnarray*}
These identities imply the full tetrahedral symmetry of the
$6j4k$-symbols, because the symmetric group on $4$ elements is
generated by the cyclic permutation $(1,2,3)$ and the transposition
$(3,4)$. To keep notations simple, we make no notational difference
between a colour weight respectively a $6j4k$-symbol and its
equivalence class.

Let $R$ be the polynomial ring over some field $\mathbb F$ whose
variables are the equivalence classes of colour weights and
$6j4k$-symbols. In this paper we will have $\mathbb F=\mathbb Q$, but
in related applications it can also be reasonable to choose for
$\mathbb F$ a finite field~\cite{kingAC}.

Let $m=|\mathcal F|$ and $n=|\mathcal G|$.  Since additional choices
play no role by the symmetry of colour weights and $6j4k$-symbols, the
following polynomial only depends on the homeomorphism type of $P$:
\[TV_{m,n}(P):=\sum_{(\varphi,\psi)\in
  \Phi_{\mathcal{F,G}}(P)} \left(\prod_{C\in \mathcal C(P)}
  w(\phi(C))\right) \cdot \left(\prod_{v\in \mathcal V(P)}
  v^{\varphi,\psi}\right) \in R\] This polynomial is the
\begriff{Turaev--Viro state sum} of $P$ of type $(m,n)$.

\section{Turaev--Viro invariants}
\label{sec:TVinvariants}

Of course, the Turaev--Viro state sum of a special spine of a compact
$3$-manifold $M$ is not yet a homeomorphism invariant of $M$, as it
will change under the $T^\pm$ moves.  Let $P_1$ be a special spine of
$M$, and let $P_2$ be obtained from $P_1$ by a single $T$ move. Let
$P_0$ be the part of $P_1$ that is unchanged by the $T$ move; we
consider $P_0$ both as a subset of $P_1$ and of $P_2$. We study how
the different terms of the state sum change under the move.

\begin{figure}[htbp]
  \centering
  {\psfrag{12}{$j$}\psfrag{13}{$j_1$}\psfrag{14}{$j_2$}
\psfrag{15}{$j_3$}\psfrag{23}{$j_4$}\psfrag{24}{$j_5$}
\psfrag{25}{$j_6$}\psfrag{34}{$j_7$}\psfrag{35}{$j_8$}
\psfrag{45}{$j_9$}
\psfrag{123}{$A_1$}\psfrag{124}{$A_2$}
\psfrag{125}{$A_3$}\psfrag{134}{$k_1$}\psfrag{135}{$k_2$}
\psfrag{145}{$k_3$}\psfrag{234}{$k_4$}\psfrag{235}{$k_5$}
 \psfrag{245}{$k_6$}\psfrag{345}{$A$}\epsfig{angle=90,file=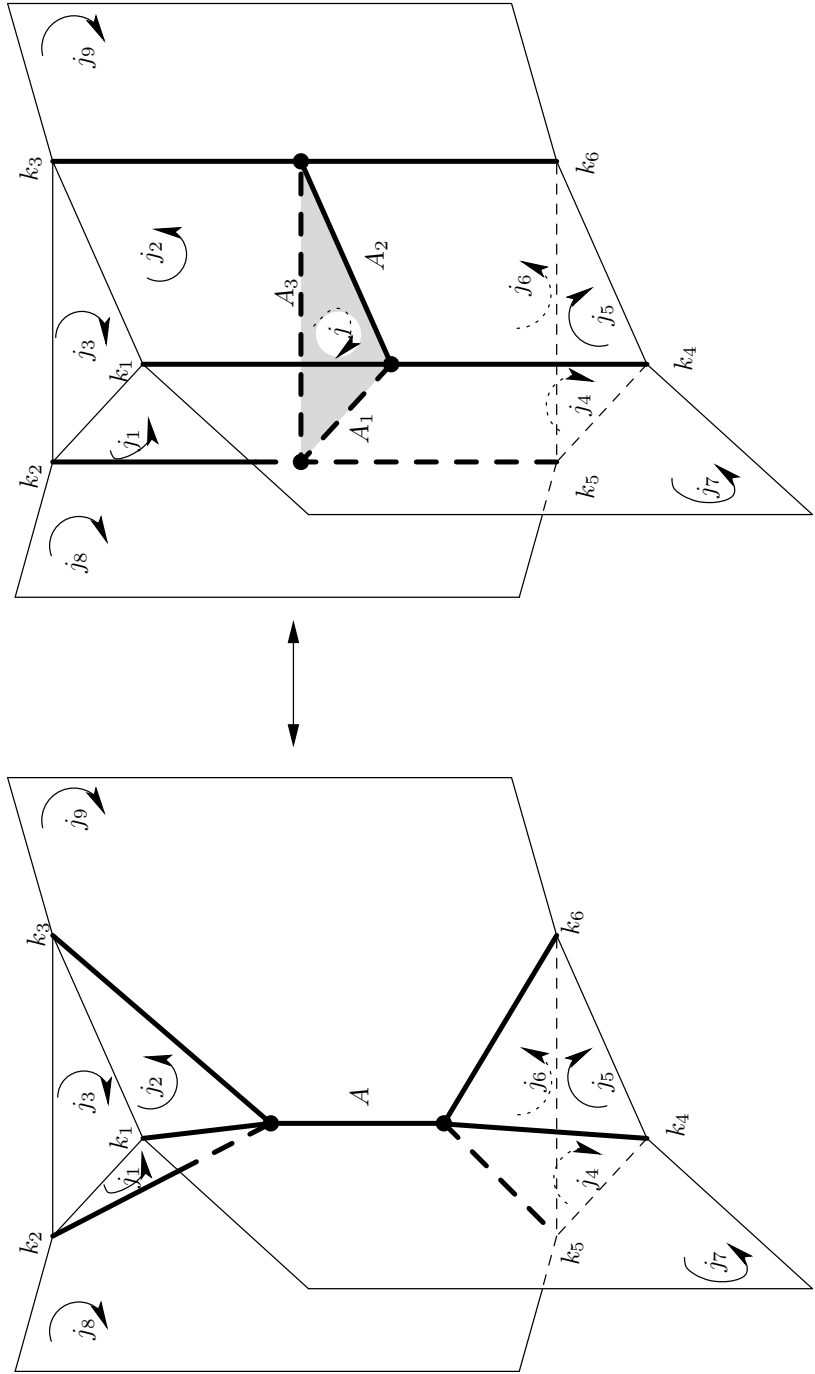}}
 \caption{Change of colourings under $T$ move}
  \label{fig:T1glg}
\end{figure}
Let $\phi_0=(\varphi_0,\psi_0)$ be the restriction of a
$\mathcal{F,G}$-coulouring $\phi_1=(\varphi_1,\psi_1)$ of $P_1$ to
$P_0$. Since any $2$-stratum of $P_1$ meets $P_0$, $\varphi_1$ is
determined by $\varphi_0$. There is one true edge of $P_1$ that is
disjoint from $P_0$, so its colour under $\psi_1$ is not determined by
$\psi_0$ --- it can be any $A\in \mathcal G$, see the left part of
Figure~\ref{fig:T1glg}. In the figure, $j_1,\dots,j_9\in \mathcal F$
denote 2-strata colours ($j_4$ and $j_6$ are the colours of the
$2$-strata that are hidden by other $2$-strata) appearing in
$\varphi_0$, and $k_1,\dots,k_6\in \mathcal G$ denote true edge
colours appearing in $\psi_0$.
If $\phi_0$ is the restriction of a $\mathcal{F,G}$-coulouring
$\phi_2=(\varphi_2,\psi_2)$ of $P_2$ to $P_0$, then $\phi_2$ is
determined by $\phi_0$, except for the colour $j\in \mathcal F$ of the
shaded triangular 2-stratum and the edge colours $A_1,A_2,A_3\in
\mathcal G$ shown in the right part of Figure~\ref{fig:T1glg}.
Let $X_{\phi_0}\in R$ be the product of the $6j4k$-symbols and colour
weights under the colouring $\phi_0$ that are associated to true
vertices and $2$-strata that are not contained in $P_1\setminus P_0$.
With an appropriate choice of orientations of $2$-strata indicated in
Figure~\ref{fig:T1glg}, we obtain
$$TV_{\mathcal{F,G}}(P_1) = \sum_{\phi_0} \sum_{A\in \mathcal G}
\jj {j_1}{j_2}{j_3}{j_7}{j_8}{j_9}^{{k_1},{k_2}}_{{k_3},{A}}\cdot \jj
{j_4}{j_5}{j_6}{-j_7}{-j_8}{-j_9}^{{k_4},{k_5}}_{{k_6},{A}}\cdot
X_{\phi_0}$$ and
{\small \begin{eqnarray*}
 &&TV_{\mathcal{F,G}}(P_2) =\\
&&\hspace*{-3mm}\sum_{\phi_0} \sum_{A_1,A_2,A_3\in
  \mathcal G}\sum_{j\in\mathcal F} w(j) \cdot \jj
{j}{j_1}{j_2}{-j_4}{-j_5}{j_7}^{{A_1},{A_2}}_{{k_1},{k_4}}\cdot \jj
{j}{j_2}{j_3}{-j_5}{-j_6}{j_9}^{{A_2},{A_3}}_{{k_3},{k_6}}\cdot \jj
{j}{j_3}{j_1}{-j_6}{-j_4}{-j_8}^{{A_3},{A_1}}_{{k_2},{k_5}}\cdot
X_{\phi_0}.\end{eqnarray*}}

Recall $m=|\mathcal F|$, $n=|\mathcal G|$.
We define the \begriff{Turaev--Viro ideal} $I_{m,n}\subset R$
of type $(m,n)$ as the ideal in $R$ that is generated by 
\begin{eqnarray*}
 &&\sum_{A\in \mathcal G}
\jj {j_1}{j_2}{j_3}{j_7}{j_8}{j_9}^{{k_1},{k_2}}_{{k_3},{A}}\cdot \jj
{j_4}{j_5}{j_6}{-j_7}{-j_8}{-j_9}^{{k_4},{k_5}}_{{k_6},{A}}
\;-\\
&& \sum_{A_1,A_2,A_3\in
  \mathcal G}\sum_{j\in\mathcal F} w(j) \cdot \jj
{j}{j_1}{j_2}{-j_4}{-j_5}{j_7}^{{A_1},{A_2}}_{{k_1},{k_4}}\cdot \jj
{j}{j_2}{j_3}{-j_5}{-j_6}{j_9}^{{A_2},{A_3}}_{{k_3},{k_6}}\cdot \jj
{j}{j_3}{j_1}{-j_6}{-j_4}{-j_8}^{{A_3},{A_1}}_{{k_2},{k_5}},
\end{eqnarray*}
for all $j_1,\dots, j_9\in \mathcal F$ and all $k_1,\dots,k_6\in
\mathcal G$.  Note that these generators are known from quantum
mechanics: If all these generators evaluate to zero, then one says
that the $6j4k$-symbols and colour weights satisfy the
``Biedenharn--Elliott equations''~\cite{landau}.
Let $tv_{m,n}(P)$ be the coset of the
Turaev--Viro state sum with respect to the Turaev--Viro ideal, i.e.
\[tv_{m,n}(P) = TV_{m,n}(P) + I_{m,n} \in R/ I_{m,n}.\] 
From
Theorem~\ref{thm:matvpierg} and the previous paragraph, we immediately
obtain
\begin{theorem}[and Definition]\label{thm:invarianz}
  If $P$ is any special spine of a compact $3$-manifold $M$ with at
  least two true vertices, then the coset $tv_{m,n}(P)$ only depends
  on the homeomorphism type of $M$. We call $tv_{m,n}(M)=tv_{m,n}(P)$
  an \begriff{ideal Turaev--Viro invariant} of $M$ of type
  $(m,n)$.
\end{theorem}

We remark that even if $M$ has a special spine $P_0$ with only one
true vertex, $tv_{m,n}(M)$ can only be computed using some special spine
with at least \emph{two} vertices. In general, $TV_{m,n}(P_0) +
I_{m,n}$ is different from $tv_{m,n}(M)$.

Note also that $tv_{m,n}(M)$ depends on the involution that we chose
for $\mathcal F$. But for simplicity we did not include the involution
in our notation. Let $N$ be the number of variables of $R$, let
$\hat {\mathbb F}$ be the algebraic closure of $\mathbb F$, and let
$\mathfrak v( I_{m,n})\subset \hat {\mathbb F}^N$ be the (affine) zero variety
associated to $I_{m,n}$. If $x\in \mathfrak v( I_{m,n})$ then, as an
obvious corollary of the preceding theorem, the state sum $TV_{m,n}(P)$
evaluated at $x$ yields an element of $\hat {\mathbb F}$ that does not
depend on the choice of a special spine $P$ with at least two vertices
of a compact $3$-manifold $M$, and we will call this a
\begriff{numerical Turaev--Viro  invariant} associated to
$tv_{m,n}(\cdot)$. By definition, if an ideal Turaev--Viro 
invariant coincides on two compact $3$-manifolds $M_1$ and $M_2$ then
\emph{all} associated numerical Turaev--Viro  invariants coincide
on $M_1$ and $M_2$.

For the following theorem, recall that the \begriff{radical} $\sqrt I$
of an ideal $I\subset R$ is the ideal formed by all polynomials $p\in
R$ with $p^n\in I$ for some $n\in\mathbb N$. An ideal is called
radical if it coincides with its radical. The zero variety of an ideal
and its radical coincide: $\mathfrak v(I)=\mathfrak v(\sqrt I)$.

\begin{definition}
  Let $M$ be a compact $3$-manifold with a special spine $P$. Let
  $tv_{m,n}(\cdot)$ be the ideal Turaev-Viro invariant obtained from
  the Turaev-Viro ideal $I_{m,n}$. The coset $$\widehat{tv}_{m,n}(M) =
  TV_{m,n}(P) + \sqrt{I_{m,n}} \in R/ \sqrt{I_{m,n}}$$
  is called the
  \begriff{universal numerical Turaev--Viro invariant} of $M$
  associated to $tv_{m,n}$.
\end{definition}

Since $I_{m,n}\subset \sqrt{I_{m,n}}$, it is clear
that $\widehat{tv}_{m,n}(M)$ is a homeomorphism invariant of
$M$. The name ``universal numerical Turaev--Viro  invariant'' is
justified by the following theorem.

\begin{theorem}\label{thm:NumVsId}
  Let $tv_{m,n}(\cdot)$ be an ideal Turaev--Viro 
  invariant.  Then for all compact $3$-manifolds $M_1$, $M_2$ holds
  $\widehat{tv}_{m,n}(M_1)=\widehat{tv}_{m,n}(M_2)$
  if and only if all numerical Turaev--Viro  invariants associated
  to $tv_{m,n}(\cdot)$ coincide on $M_1$ and $M_2$.
\end{theorem}
\begin{proof}
  If
  $\widehat{tv}_{m,n}(M_1)=\widehat{tv}_{m,n}(M_2)$
  then all numerical Turaev--Viro  invariants associated to
  $tv_{m,n}(\cdot)$ coincide on $M_1$ and $M_2$, since
  $\mathfrak v(I_{m,n}) = \mathfrak
  v(\sqrt{I_{m,n}})$.
  
  Now assume that for special spines $P_1$, $P_2$ of $M_1$, $M_2$ with
  at least two true vertices holds
  $TV_{m,n}(P_1)(x)=TV_{m,n}(P_2)(x)$ for all
  $x\in\mathfrak v(I_{m,n})$. So the polynomial
  $\left(TV_{m,n}(P_1)-TV_{m,n}(P_2)\right)\in
  R$
  vanishes on $\mathfrak v(I_{m,n})$. Hence by Hilbert's
  Nullstellensatz (in the formulation stated
  in~\cite[Corollary~2.6.17]{KreuzerRobbiano}) we have 
  $$
  \left(TV_{m,n}(P_1)-TV_{m,n}(P_2)\right)\in
  \sqrt{I_{m,n}},$$ and therefore
  $\widehat{tv}_{m,n}(P_1)=\widehat{tv}_{m,n}(P_2)$.
\end{proof}

Actually it turns out that the Turaev--Viro ideals studied in
Section~\ref{sec:expl} and~\ref{sec:computations} are not radical.  So
we have reason to expect that, in general, with an ideal Turaev--Viro
invariant one can distinguish strictly more manifolds than with all
its associated numerical Turaev--Viro invariants together.

It is well known that numerical Turaev--Viro  invariants exist for
arbitrarily large colour sets. First examples have been presented by
V.~Turaev and O.~Viro~\cite{TuraevViro92}. The representation theory
of quantum groups yields a very successful machinery for constructing
numerical Turaev--Viro  invariants.

We come to another potential application of ideal Turaev--Viro 
invariants. Let $\tilde c(M)$ be the minimal number of true vertices
of a special spine of a compact $3$-manifold $M$. This is related to
Matveev's notion of \begriff{complexity} of manifolds, $c(M)$: If $M$
is a closed irreducible $3$-manifold different from the $3$-sphere,
the projective space and the lens space $L(3,1)$ then $c(M)=\tilde
c(M)$.

Let $p\in R$ be a polynomial. Let $\deg_w(p)$ be the total degree of
$p$ in the colour weights, and let $\deg_{6j}(p)$ be the total degree
of $p$ in the $6j4k$-symbols. For any subset $A\subset R$, let
$\deg_w(A)=\min\{\deg_w(p)\co p\in A\}$ and
$\deg_{6j}(A)=\min\{\deg_{6j}w(p)\co p\in A\}$.

\begin{lemma}\label{lem:untereschranke}
  Let $tv_{m,n}(\cdot)$ be an ideal Turaev--Viro  invariant. For
  any closed $3$-manifold $M$ with $\tilde c(M)>1$, we have
  $$ \tilde c(M)\ge \max\left\{\deg_w\left(tv_{m,n}(M)\right) -1,
                                    \deg_{6j}\left(tv_{m,n}(M)\right)\right\}.$$
\end{lemma}
\begin{proof}
  Let $P$ be a special spine of $M$ with $\tilde c(M)$ true vertices.
  Since $\tilde c(M)>1$, $P$ has at least two true vertices, hence we
  can compute $tv_{m,n}(M)$ using $P$. Since $\d M=\emptyset$, we have
  $|\mathcal C(P)| =\tilde c(M)+1$. Hence for the polynomial
  $TV_{m,n}(P)\in tv_{m,n}(M)$ we find
  $\deg_w\left(TV_{m,n}(P)\right)\le \tilde c(M)+1$ and
  $\deg_{6j}\left(TV_{m,n}(P)\right)\le \tilde c(M)$ (possibly with
  strict inequality if simplifying assumption apply for
  $tv_{m,n}(\cdot)$; see Section~\ref{sec:assumptions}).
\end{proof}

\section{Computation of ideal Turaev--Viro  invariants}
\label{sec:implementation}

How can one distinguish manifolds using ideal Turaev--Viro 
invariants? The first task is to present a special $2$-polyhedron in a
form that is accessible for computers.
S.~Matveev~\cite[Sec.~7.1]{matveev} introduced a way of encoding
special $2$-polyhedra by lists of cyclic sequences of integers. This
is roughly as follows. Let $P$ be a special $2$-polyhedron. We number
the true edges of $P$ and provide them with an arbitrary orientation.
Let $C$ be an oriented $2$-stratum of $P$. When we track the oriented
boundary of $C$, we obtain a cyclic sequence of the oriented true
edges met by $\d C$.  This cyclic sequence of oriented edges is
encoded by a cyclic sequence of integers whose absolute value gives
the number of each edge met by $\d C$, with positive (resp. negative)
sign if the orientation of $\d C$ and of the edge coincides (resp.
does not coincide). It turns out that this list of cyclic sequences of
integers determines $P$ up to homeomorphism.

It is not difficult to deduce from the representation of $P$ how the
oriented $2$-strata and true edges of $P$ meet at the true vertices.
Hence, one can easily implement the computation of the Turaev--Viro
state sum. If $P$ is a special spine of a closed $3$-manifold and has
$c$ true vertices, then $|C(P)|=c+1$ and $|E(P)|=2c$. Again, let
$m=|\mathcal F|$ and $n=|\mathcal G|$. The number of summands in
$TV_{m,n}(P)$ is roughly $m^{c+1}\cdot n^{2c}$. So the computation of
the state sum is easy but for large $c$ quite time-consuming.

Let $P_1$ and $P_2$ be special spines of compact $3$-manifolds $M_1$
and $M_2$. We are now able to compute $TV_{m,n}(P_1)$ and
$TV_{m,n}(P_2)$. But how can we determine whether
$tv_{m,n}(M_1)=tv_{m,n}(M_2)$ or not? In other words, we need to
compare cosets with respect to ideals in a polynomial ring over a
field. This is algorithmically possible by the theory of Gröbner
bases. For an introduction to that subject, we refer the reader
to~\cite{Froberg} or~\cite{KreuzerRobbiano}, among many other possible
sources.

Firstly, we need to choose an \begriff{admissible monomial ordering}
$<$ on $R$; this is a total order on the set of monomials (i.e.,
products of variables) of $R$ such that $1<\mathfrak m$ for any
monomial $\mathfrak m\in R$ and such that $\mathfrak m_1<\mathfrak
m_2$ implies $\mathfrak m\mathfrak m_1<\mathfrak m\mathfrak m_2$ for
all monomials $\mathfrak m,\mathfrak m_1,\mathfrak m_2\in R$. For a
polynomial $f\in R$, the \begriff{leading monomial} of $f$ with
respect to $>$ is denoted by $lm_>(f)$. If an admissible ordering on
$R$ is given then one can generalise the usual division algorithm of univariate
polynomials to multivariate polynomials and can define the remainder
$\mathrm{rem}_>(f,g)\in R$ of a polynomial $f\in R$ with respect to a
polynomial $g\in R$. In general the remainder will depend on the
chosen ordering.

Let $I=\langle g_1,\dots,g_k\rangle\subset R$ be the ideal generated
by the polynomials $g_1,\dots, g_k\in R$. If one wants to test whether
some polynomial $f\in R$ belongs to $I$, it is a reasonable idea to
iteratively compute the remainder $\mathrm{rem}_>(f;g_1,\dots,g_k)$ of
$f$ with respect to $g_1,\dots,g_k$, i.e., $\mathrm{rem}_>\left(\dots
  \mathrm{rem}_>\left(\mathrm{rem}_>(f,g_1),g_2\right)\dots,g_k\right)$.
Certainly if $\mathrm{rem}_>(f;g_1,\dots,g_k)=0$ then $f\in I$.
However, in general the converse is not true. Moreover, in general
$\mathrm{rem}_>(f;g_1,\dots,g_k)$ depends on the order of $g_1,\dots,
g_k$.

A \begriff{Gröbner base} of $I$ with respect to $>$ is a finite set
$B\subset I$ such that $$\left\langle \{lm(f)\co f\in B\}\right\rangle =
\left\langle \{lm(f)\co f\in I\}\right\rangle.$$ It turns out that any
Gröbner base of $I$ is a generating subset of $I$, and that any ideal
in $R$ \emph{has} a Gröbner base (specifically, any ideal is finitely
generated). If $B$ satisfies some additional hypothesis
(see~\cite[Sec.~3.7]{Froberg} for details), it is called \emph{reduced} Gröbner
base, and turns out to be unique, hence depends only on $I$ and $>$.
The reduced Gröbner base can be algorithmically constructed, given
{an} arbitrary finite generating subset of $I$. One of the main
features of (not necessarily reduced) Gröbner bases is that they allow
the computation of a unique representative for any coset $f+I\in R/I$:
If $\{b_1,\dots,b_n\}$ is a Gröbner base of $I$ with respect to $>$,
then for any $f,g\in R$ one has $\mathrm{rem}_>(f;b_1,\dots,b_n)=
\mathrm{rem}_>(g;b_1,\dots,b_n)$ if and only if $f+I=g+I$. Moreover,
$\mathrm{Nf}_>(f,I)=\mathrm{rem}_>(f;b_1,\dots,b_n)$ does not depend on the
choice of a Gröbner base for $I$ or on the order of $b_1,\dots,b_n$,
but only on $f+I$ and $>$, and is therefore called the \begriff{normal
  form} of $f+I$ with respect to $>$.
The computation of Gröbner bases and normal forms is implemented in
various computer algebra systems, e.g. in \texttt{Maple V} or
\texttt{Mathematica}, and there is also specialised software like
\texttt{Singular}, \texttt{Macaulay~2} or \texttt{bergman}.

Our computations involve the following three steps.
\begin{enumerate}
\item Produce the list of variables of $R$ and the list of generators of
  $I_{m,n}$ defined in the previous section.
\item Compute a Gröbner base of $I_{m,n}$ for the
  chosen admissible monomial ordering $>$.
\item For any special $2$-polyhedron $P$, compute
  $TV_{m,n}(P)$, and 
\item compute the normal form of $TV_{m,n}(P)+I_{m,n}$ using the
  Gröbner base obtained in step~2.
\end{enumerate}
For step~1 and~3, we wrote \texttt{maple~V} programs. For step~2
and~4, we used \texttt{Singular}. If we want to compute the universal
numerical Turaev--Viro invariant associated to $tv_{m,n}$, we simply
replace $I_{m,n}$ by $\sqrt{I_{m,n}}$, which is possible since one can
compute a finite set of generators of $\sqrt{I_{m,n}}$ for any finite
set of generators of $I_{m,n}$ (we used \texttt{Singular} for that
purpose).

We conclude this section with some remarks on the computational
complexity. The number of $6j4k$-symbols is $m^6\cdot
n^4$, and we have $m$ colour weights. The number
of variables of $R$ is slightly less since we take equivalence
classes, though it still grows rapidly with the number of colours. We
have roughly $m^9\cdot n^6$ generators of
$I_{m,n}$ (some of them coincide by symmetry), which
are polynomials of degree $4$ involving up to $n +
m \cdot n^3$ monomials. Apart from the sheer
size of the generating system of $I_{m,n}$, step~1
is easy. And so is step~3, except for the size of the state sum and
for the number of different manifolds that we want to compute
invariants for.

Step~2 is the most critical one. The main problem for computing Gröbner
bases is the number of variables, especially if a lexicographic order
is used. So one can not expect to be able to use large colour sets.
However, the examples exposed in Section~\ref{sec:expl} show that
to some extent the computation of ideal Turaev--Viro  invariants is
feasible. Fortunately the most critical step~2 is to be performed only
once, it is not needed to repeat it for any manifold.

\section{Simplifying assumptions}
\label{sec:assumptions}

One way to overcome the complexity problems mentioned in the previous
section is to introduce simplifying assumptions. For instance, we can
restrict the set of colourings by sending some of the $6j4k$-symbols
to zero. This strategy is supported by the numerical Turaev--Viro 
invariants obtained from quantum groups. E.g., in the invariants
constructed in~\cite{TuraevViro92}, we have colour sets $\mathcal
F=\{0,\dots,k-1\}$ and $\mathcal G=\{*\}$, and the $6j$-symbols vanish
unless around any true edge the three 2-strata colours satisfy
triangle inequalities and have even sum (``admissible'' colouring). In
the numerical examples, usually one has some colour $z\in \mathcal F$
such that $w(z)=\jj zzzzzz=1$, and we could assume the same. We even
could adopt the values of some more $6j4k$-symbols occuring in
numerical Turaev--Viro  invariants, and only declare the remaining
symbols as variables. Since the number of variables is critical for
computing Gröbner bases, this is an efficient strategy.

\begin{lemma}\label{lem:scaling}
  If $\jj abcdef$ are the $6j$-symbols of a numerical Turaev--Viro 
  invariant of type $(m,1)$ with $\mathcal G=\{*\}$ then one obtains a
  numerical Turaev--Viro  invariant of type $(m,n)$ by maintaining
  the colour weights and defining the $6j4k$-symbols by scaling the
  $6j$-symbols: $\jj abcdef^{A,B}_{C,D} = \frac 1{n^2}\cdot\jj abcdef$
  for all $A,B,C,D\in \mathcal G$.
\end{lemma}
\begin{proof}
  The scaling implies $$\sum_{A\in \mathcal G} \jj
{j_1}{j_2}{j_3}{j_7}{j_8}{j_9}^{{k_1},{k_2}}_{{k_3},{A}}\cdot \jj
{j_4}{j_5}{j_6}{-j_7}{-j_8}{-j_9}^{{k_4},{k_5}}_{{k_6},{A}} = \frac n{n^4}
\jj {j_1}{j_2}{j_3}{j_7}{j_8}{j_9}\cdot \jj
{j_4}{j_5}{j_6}{-j_7}{-j_8}{-j_9}$$
and \begin{align*}
\sum_{A_1,A_2,A_3\in
  \mathcal G}&\sum_{j\in\mathcal F} w(j) \cdot \jj
{j}{j_1}{j_2}{-j_4}{-j_5}{j_7}^{{A_1},{A_2}}_{{k_1},{k_4}}\cdot \jj
{j}{j_2}{j_3}{-j_5}{-j_6}{j_9}^{{A_2},{A_3}}_{{k_3},{k_6}}\cdot \jj
{j}{j_3}{j_1}{-j_6}{-j_4}{-j_8}^{{A_3},{A_1}}_{{k_2},{k_5}}  \\
=&\frac {n^3}{n^6}\sum_{j\in\mathcal F} w(j) \cdot \jj
{j}{j_1}{j_2}{-j_4}{-j_5}{j_7}\cdot \jj
{j}{j_2}{j_3}{-j_5}{-j_6}{j_9}\cdot \jj
{j}{j_3}{j_1}{-j_6}{-j_4}{-j_8}.
\end{align*}
So the generators of $I_{m,n}$ are obtained from the
generators of $I_{m,1}$ by scaling with $\frac 1{n^3}$.
Hence if the colour weights and the $6j$-symbols correspond to a point
of $\mathfrak v(I_{m,1})$ then the colour weights and the
$6j4k$-symbols correspond to a point of $\mathfrak v(I_{m,n})$.
\end{proof}
Of course, the preceding Lemma will not yield an essentially new
invariant. But it suggests to give \emph{some} of the $6j4k$-symbols
the value $\jj abcdef^{A,B}_{C,D} = \frac 1{n^2}\cdot\jj abcdef$ and
keep the remaining symbols as variables.

Finally, we can reduce the number of variables of $R$ by providing the
$6j4k$-symbols with additional symmetries. For instance, we could
assume 
\begin{equation}\label{eqn:symmetry} 
  \jj abcdef^{A,B}_{C,D} = \jj abcdef^{\sigma(A),\sigma(B)}_{\sigma(C),\sigma(D)},
\end{equation}
for any $a,b,c,d,e,f\in \mathcal F$, any $A,B,C,D\in \mathcal G$ and
any permutation $\sigma$ of $\{A,B,C,D\}$.

The general machinery always remains the same: We have a polynomial
ring whose variables correspond to $6j4k$-symbols and colour weights,
we have a state sum associated to any special $2$-polyhedron, and we
have an ideal such that the coset of the state sum does not change
under $T^\pm$ moves.  Therefore we still call the resulting
homeomorphism invariant of compact $3$-manifolds ``an ideal
Turaev--Viro  invariant $\widetilde{tv}_{m,n}(\cdot)$ of type
$(m,n)$'' (we mark by the tilde the existence of simplifying
assumptions), and it should be clear how to define the notion of a
(universal) numerical Turaev--Viro  invariant associated to
$\widetilde{tv}_{m,n}(\cdot)$ --- Theorem~\ref{thm:NumVsId} holds with
the obvious changes also in the new setting. 

Because all mentioned simplifying assumptions are known to hold for
some numerical Turaev--Viro  invariants, we can be sure that
\emph{ideal} Turaev--Viro  invariants subject to these assumptions
are nontrivial as well.
%

One may also think of adding \emph{more} generators to $I_{m,n}$. This
might accelerate the computation of a Gröbner base, because additional
polynomials can be used to simplify the generators of $I_{m,n}$.  A
source of additional generators could be an adaption of a Gröbner base
of a Turaev--Viro ideal of type $(m,1)$.

\section{Examples}
\label{sec:expl}

In this section, we present several ideal Turaev--Viro invariants.  Our
computational results on these invariants will be stated in the next
section. For all examples, we chose $\mathbb F=\mathbb Q$.

Our first example, $\widetilde{tv}_{2,1}$, is constructed similarly to
Matveev's
$\eps$-invariant~\cite{MatveevOvchinnikovSokolov},~\cite{matveev}.
Actually the $\eps$-invariant is obtained as evaluation of
$\widetilde{tv}_{2,1}$ (see below). For $\widetilde{tv}_{2,1}$, we have $\mathbb
F=\mathbb Q$, $2$-strata colours $\mathcal F=\{1,2\}$ and trivial edge
colours $\mathcal G=\{*\}$. So we have $6j$-symbols rather than
$6j4k$-symbols. We work under the simplifying assumption that the
$6j$-symbol of a vertex $v$ vanishes if there is some true edge $e$
meeting $v$ so that there are exactly two $2$-strata of colour $1$
meeting $e$ (counted with multiplicity).  We also assume $w(1)=1$ and
$\jj 111111=1$, which holds for the $\eps$-invariant as well. So the
variables of $R$ are four equivalence classes of $6j$-symbols and the
colour weight $w(2)$, and we provide $R$ with degree reverse
lexicographic order, where $\jj 112122> \jj 212212> \jj 212222> \jj
222222>w(2)$.

The Turaev--Viro ideal in this setting is generated by 12 polynomials:
{\allowdisplaybreaks\begin{align*} & \jj 222222^2- \jj
    212222^3-w(2)\cdot \jj 222222^3, \\ & \jj 212222^2- \jj
    212222^2\cdot \jj 212212-w(2)\cdot \jj 222222^2\cdot \jj 212222,
    \\ & \jj 212222\cdot \jj 212212-w(2)\cdot \jj 212222^3, \\ & \jj
    212212\cdot \jj 112122-w(2)\cdot \jj 212212^2\cdot \jj 112122, \\
    & \jj 212222\cdot \jj 112122-w(2)\cdot \jj 212222^2\cdot \jj
    112122, \\ & \jj 212212^2- \jj 212212\cdot \jj 112122^2, \\ & \jj
    112122^2- \jj 212212^3-w(2)\cdot \jj 212222^3, \\ & \jj
    112122-w(2)\cdot \jj 112122^3, \\ & -\jj 212222\cdot \jj
    212212^2-w(2)\cdot \jj 222222\cdot \jj 212222^2, \\ & \jj
    222222\cdot \jj 212222-w(2)\cdot \jj 222222\cdot \jj 212222^2, \\
    & \jj 212222^2- \jj 212222\cdot \jj 112122^2, \\ & \jj
    212222^2-w(2)\cdot \jj 212222^2\cdot \jj 212212,
\end{align*}}
compare~\cite[Sec.~8.1.2]{matveev}. Note that without the assumption
$\jj 111111=w(1)=1$ one has two additional generators.  One obtains a
(non-reduced) Gröbner base formed by 22 polynomials:
{\allowdisplaybreaks\begin{align*}
&   \jj 212222\cdot  \jj 212212-
 \jj 212222^2, \\ & - \jj 112122+
 \jj 112122\cdot w(2)^2+ \jj 212212\cdot  \jj 112122-
 \jj 212222\cdot  \jj 112122-w(2)\cdot  \jj 112122, \\ & 
 \jj 212222^2+w(2)\cdot  \jj 222222\cdot  \jj 212222, \\ &
 w(2)\cdot  \jj 212222^2- \jj 212222^2+
 \jj 222222\cdot  \jj 212222, \\ &  \jj 212212\cdot 
 \jj 112122+w(2)\cdot  \jj 212222\cdot 
 \jj 112122-w(2)\cdot  \jj 112122, \\ & w(2)\cdot 
 \jj 212212^2- \jj 112122^2, \\ & w(2)\cdot 
 \jj 212212\cdot  \jj 112122- \jj 112122, \\ & w(2)\cdot 
 \jj 112122^2-w(2)\cdot  \jj 212222^2-
 \jj 212212^2, \\ & \jj 222222^2\cdot  \jj 212222+
 \jj 212222^2+2\cdot  \jj 222222\cdot  \jj 212222, \\ & 
 \jj 212222^2+ \jj 222222\cdot  \jj 212222^2+
 \jj 222222\cdot  \jj 212222, \\ &  \jj 112122\cdot 
 \jj 222222\cdot  \jj 212222+ \jj 212212\cdot 
 \jj 112122+2\cdot  \jj 212222\cdot  \jj 112122-w(2)\cdot 
 \jj 112122, \\ &  \jj 222222\cdot  \jj 212212\cdot 
 \jj 112122- \jj 112122\cdot  \jj 222222\cdot 
 \jj 212222\\&
 \mbox{\hspace*{3cm}}-w(2)\cdot  \jj 112122\cdot  \jj 222222-
 \jj 212222\cdot  \jj 112122, \\ & 
 \jj 212222^3-w(2)\cdot  \jj 212222^2+
 \jj 212222^2, \\ &  \jj 112122\cdot  \jj 212222^2-
 \jj 112122\cdot  \jj 222222\cdot  \jj 212222-
 \jj 212222\cdot  \jj 112122, \\ & - \jj 212222^2+
 \jj 212222\cdot  \jj 112122^2, \\ &  \jj 212212^3-
 \jj 112122^2+ \jj 212222\cdot  \jj 212212, \\ & 
 \jj 112122\cdot  \jj 212212^2+ \jj 212222\cdot 
 \jj 112122- \jj 112122, \\ & - \jj 212212^2+
 \jj 212212\cdot  \jj 112122^2, \\ & 
 \jj 112122^3+w(2)\cdot  \jj 212222\cdot 
 \jj 112122-w(2)\cdot  \jj 112122, \\ & w(2)\cdot 
 \jj 222222^3+ \jj 212222^3- \jj 222222^2, \\ & 
 \jj 222222^3\cdot  \jj 112122-w(2)\cdot  \jj 112122\cdot 
 \jj 222222^2-4\cdot  \jj 212212\cdot  \jj 112122
\\ &\qquad\qquad \mbox{}-6\cdot 
 \jj 212222\cdot  \jj 112122+4\cdot w(2)\cdot 
 \jj 112122, \\ &  \jj 222222^3\cdot  \jj 212212^2-
 \jj 112122^2\cdot  \jj 222222^2+4\cdot 
 \jj 212222^2+6\cdot  \jj 222222\cdot  \jj 212222.
\end{align*}}

According to eqn.~(8.5) in~\cite{matveev}, one obtains the
$\eps$-invariant by evaluation as follows:
\[
  \jj 112122 = \eps^{-\frac 12},
  \jj 212212 = \eps^{-1},
 \jj 212222 = \eps^{-1}, 
\jj 222222 = -\eps^{-2},
 w(2)=\eps,
\]
where $\eps$ is any root of $\eps^2-\eps-1$.

This is the only example for which we write down the defining
polynomials in this article. In all other cases, there are simply too
many polynomials, and we refer to the material that we provide on our
web site~\cite{web}. Note that one obtains a smaller Gröbner base with
a lexicographic order on $R$. However we did not succeed to work with
lexicographic orders in our other examples.

By a result of Matveev--Nowik~\cite{MatveevNowik}, there are pairs of
non-homeomorphic compact manifolds that can not be distinguished by
any Turaev--Viro  invariant of type $(m,1)$ with a trivial
involution on the set of $2$-strata colours (i.e., $f=-f$ for all
$f\in \mathcal F$). The result was formulated for numerical
invariants, but it readily applies for ideal invariants as well.
Therefore we also constructed two invariants of type $(3,1)$ with
non-trivial involution and an invariant of type $(2,2)$.

For type $(3,1)$, let $\mathcal F=\{-1,0,1\}$ with the usual
involution $-(-1)=1$, $-0=0$. We assume for simplification that
$w(0)=1$ and $\jj 000000=1$. Then, we have 41 equivalence classes of
$6j$-symbols and one remaining colour weight $w(1)$. We obtain $1661$
generators for the Turaev--Viro ideal, and after two days
\texttt{Singular} succeeds with finding a Gröbner base with respect to
some degree reverse lexicographic order formed by $1297$ polynomials. We
denote the resulting ideal Turaev--Viro  invariant by $tv_{3,1}^+(\cdot)$,
where the ``$+$'' shall denote that the involution on $\mathcal F$ is
non-trivial.

In order to compare an ideal invariant of a given type with another
ideal invariant of the same type subject to simplifying assumptions,
we also consider the following setting. Again, let $\mathcal
F=\{-1,0,1\}$ with the usual involution $-(-1)=1$, $-0=0$. We assume
$w(0)=1$, but we do not assume $\jj 000000=1$. Instead, we assume
that the $6j$-symbol of a vertex $v$ vanishes if there is some true
edge $e$ meeting $v$ so that there are exactly two $2$-strata of
colour $0$ meeting $e$ (counted with multiplicity). Hence, the
assumption is essentially the same as in the case of
$\widetilde{tv}_{2,1}$. There remain $21$ equivalence classes of
$6j$-symbols. The Turaev--Viro ideal is generated by $474$ polynomials,
and after a few seconds \texttt{Singular} finds a Gröbner base of
$337$ polynomials. We denote the resulting invariant by
$\widetilde{tv}_{3,1}^+(\cdot)$.

We now come to an ideal Turaev--Viro invariant of type $(2,2)$. We use
$\mathcal F=\{1,2\}$ with trivial involution, and $\mathcal
G=\{1,2\}$. For simplification, we assume $w(1)=\jj
111111^{A,B}_{C,D}=\frac 14$, which is justified by
Lemma~\ref{lem:scaling}, since in $\widetilde{tv}_{2,1}(\cdot)$ we
have $\jj 111111=1$. 
We assume that the $6j$-symbol of a vertex $v$ 
vanishes if there is some true edge $e$ meeting $v$ so that there are
exactly two $2$-strata of colour $1$ meeting $e$ (counted with
multiplicity).  
Moreover, we assume the additional symmetry stated in
Equation~(\ref{eqn:symmetry}). We then have $22$ equivalence classes of
$6j4k$-symbols. The Turaev--Viro ideal is generated by $353$
polynomials, but we did not succeed to compute a Gröbner base in this
setting. Therefore we enlarged the ring $R$ by a new variable $X$ and
enlarged the Turaev--Viro ideal by adding $22$ generators, obtained as
follows: We define $\jj abcdef=X\cdot \jj abcdef^{1,1}_{1,1}$ for
$a,\dots,f\in \mathcal F$, apply this definition to the $6j$-symbols
in the $22$ polynomials in the above Gröbner base used to compute
$\widetilde{tv}_{2,1}$, and append them to the list of generators of
the Turaev--Viro ideal.
After this enlargement of the ideal, \texttt{Singular} finds a Gröbner
base formed by $449$ polynomials within a few minutes. We denote the
resulting invariant by $\widetilde{tv}_{2,2}(\cdot)$.

\section{Computational results}
\label{sec:computations}

In this section, we report the results of computing the invariants
presented in the previous section on lists of closed orientable
manifolds.  We computed $\widetilde {tv}_{2,1}(\cdot)$, $\widetilde
{tv}_{3,1}^+(\cdot)$ and $tv_{3,1}^+(\cdot)$ for irreducible manifolds
up to complexity $9$, and $\widetilde{tv}_{2,2}(\cdot)$ up to
complexity $6$. The following list of statements is result of our
computations.
\begin{proposition}\label{thm:results}
  \mbox{}\begin{enumerate}
  \item The Turaev--Viro ideals involved in the construction of
    $\widetilde {tv}_{2,1}(\cdot)$, $\widetilde {tv}_{3,1}^+(\cdot)$ and 
    $tv_{3,1}^+(\cdot)$ 
    are not radical.
  \item On the $1900$ closed irreducible orientable $3$-manifolds of complexity
    $\le 9$, the $\eps$-invariant takes $35$ different values,
    $\widetilde {tv}_{2,1}(\cdot)$ takes $134$ different values and
    $\widetilde {tv}_{3,1}^+(\cdot)$ takes $242$ different values,
    whereas homology takes $272$ different values.
  \item Using the combination of homology and $\widetilde
    {tv}_{3,1}(\cdot)$ one can distinguish $764$ homeomorphism types
    of closed irreducible orientable $3$-manifolds of complexity $\le 9$.
  \item On closed irreducible orientable  $3$-manifolds of complexity $\le 9$,
    $\widetilde {tv}_{3,1}^+(\cdot)$ and $tv_{3,1}^+(\cdot)$ are
    equivalent invariants. On closed irreducible orientable $3$-manifolds of
    complexity $\le 6$, $\widetilde {tv}_{2,1}(\cdot)$ and
    $\widetilde{tv}_{2,2}(\cdot)$ are equivalent invariants.
  \item On the closed irreducible orientable $3$-manifolds that we
    considered, the ideal Turaev--Viro invariants $\widetilde
    {tv}_{2,1}(\cdot)$, $\widetilde {tv}_{3,1}^+(\cdot)$ and
    $tv_{3,1}^+(\cdot)$ are equivalent to their associated universal
    numerical Turaev--Viro invariant.
  \item The lower bound for the complexity stated
    in~\ref{lem:untereschranke} is trivial in all examples that we
    computed.
  \item Ideal Turaev--Viro invariants are, in general, not
    multiplicative under connected sum of compact $3$-manifolds.
  \end{enumerate}
\end{proposition}
\begin{proof}
  Statements (1)--(6) are simply obtained by running the algorithm
  scetched in Section~\ref{sec:implementation} on lists of special
  spines of manifolds. Statement~(6) holds since in fact the normal
  forms of the invariants have degree $2$ in the $6j4k$-symbols and
  degree $3$ in the colour weights. 
  
  We go in a little more detail with Statement~(1), in the case of
  $\widetilde{tv}_{2,1}(\cdot)$. A Gröbner base of $\sqrt{\tilde I_{2,1}}$ is given by 
 {\allowdisplaybreaks \begin{align*}
    &\jj 222222 \cdot \jj 212222 -\jj 212222\cdot w(2)+ 2\cdot \jj
    212222, \\ & \jj 212222^2- \jj 212222 \cdot w(2)+\jj 212222, \\ &
    \jj 212222 \cdot \jj 212212 -\jj 212222 \cdot w(2)+\jj 212222, \\ &
    \jj 212212^2- \jj 212212\cdot w(2)+\jj 212222, \\ & \jj 112122^2-
    \jj 212212, \\ & \jj 212222 \cdot w(2)^2-\jj 212222 \cdot w(2)-\jj
    212222, \\ & \jj 212212\cdot w(2)^2- \jj 212222 \cdot w(2)-\jj
    212212, \\ & -\jj 112122 +\jj 112122 \cdot w(2)^2+\jj 212212\cdot
    \jj 112122 -\jj 212222 \cdot \jj 112122 -w(2)\cdot \jj 112122, \\ &
    \jj 222222^2\cdot w(2)-\jj 212222 -\jj 222222, \\ & \jj
    212212\cdot \jj 112122 +w(2)\cdot \jj 212222 \cdot \jj 112122
    -w(2)\cdot \jj 112122, \\ & w(2)\cdot \jj 212212\cdot \jj 112122 -
    \jj 112122 \cdot w(2)^2-\jj 212212\cdot \jj 112122\\&\qquad\qquad + \jj 212222
    \cdot \jj 112122 +w(2)\cdot \jj 112122, \\ & \jj 212212\cdot \jj
    222222^2-\jj 212212\cdot \jj 222222 \cdot w(2)+2\cdot \jj 212222
    \cdot w(2)- 4\cdot \jj 212222, \\ & \jj 112122 \cdot \jj 222222^2-
    w(2)\cdot \jj 112122 \cdot \jj 222222 +\jj 112122 \cdot
    w(2)^2\\&\qquad\qquad +3\cdot \jj 212212\cdot \jj 112122 +\jj 212222 \cdot \jj
    112122 -3\cdot w(2)\cdot \jj 112122 -\jj 112122, \\ & \jj 222222
    \cdot \jj 212212\cdot \jj 112122 - w(2)\cdot \jj 112122 \cdot \jj
    222222 +\jj 212212\cdot \jj 112122\\&\qquad\qquad +\jj 212222 \cdot \jj 112122
    -w(2)\cdot \jj 112122,
  \end{align*}}
and when one computes the normal form of these generators using the above mentioned
Gröbner base for $\tilde I_{2,1}$, one sees that $\sqrt{\tilde
  I_{2,1}}\not\subset \tilde I_{2,1}$, since the following
non-vanishing normal forms remain: 
{\allowdisplaybreaks \begin{align*}
    &\jj 212212\cdot \jj 222222^2-\jj 212212\cdot \jj 222222\cdot
    w(2)+2\cdot \jj 212222\cdot w(2) -4\cdot \jj 212222,\\ & \jj
    112122^2-\jj 212212, \\ & \jj
    222222\cdot \jj 212222-\jj 212222\cdot w(2)+ 2\cdot \jj 212222, \\
    & \jj 212222^2-\jj 212222\cdot w(2)+\jj 212222, \\ & \jj
    212212^2-\jj 212212 \cdot w(2)+\jj 212222, \\ &\jj 212222\cdot
    w(2)^2- \jj 212222\cdot w(2)-\jj 212222, \\ & \jj 212212\cdot
    w(2)^2-\jj 212222\cdot w(2)-\jj 212212, \\ & \jj 222222^2\cdot
    w(2)-\jj 212222-\jj 222222, \\ & \jj 112122\cdot \jj
    222222^2-w(2)\cdot \jj 112122\cdot \jj 222222+2\cdot \jj
    212212\cdot \jj 112122\\ &\qquad\qquad +2\cdot \jj 212222\cdot \jj
    112122-2\cdot w(2)\cdot \jj 112122.
  \end{align*}}
  
  We show statement~(7) again for $\widetilde{tv}_{2,1}(\cdot)$, but
  it holds analogously also for our other examples of ideal
  Turaev--Viro invariants. We found
\begin{align*}
    \mathrm{Nf}_>&\left(\widetilde{tv}_{2,1}\left(L(7,2)\right)\right)\\ &= w(2)^3\cdot\jj 222222^2 -\jj 212222^2
                        + \jj 222222\cdot\jj 212222+ 1\\
     \mathrm{Nf}_>&\left(\widetilde{tv}_{2,1}\left(L(8,3)\right)\right)\\ &= w(2)^3\cdot\jj 222222^2+2\cdot\jj 212222^2
                          -\jj 222222\cdot\jj 212222+1,
\end{align*}    
   but
\begin{align*}
     \mathrm{Nf}_>&\left(\widetilde{tv}_{2,1}\left(L(8,3)\#L(8,3)\right)\right)\\ &=w(2)^3\cdot\jj 222222^2+3\cdot\jj 212212^2
                                     +7\cdot\jj 212222^2
                             -6\cdot\jj 222222\cdot\jj 212222+1\\
  \not= \mathrm{Nf}_>&\left(\mathrm{Nf}_>\left(\widetilde{tv}_{2,1}\left(L(8,3)\right)\right)
               \cdot \mathrm{Nf}_>\left(\widetilde{tv}_{2,1}\left(L(8,3)\right)\right)\right) \\
  &= w(2)^4\cdot\jj 222222^2+2\cdot w(2)^3\cdot\jj 222222^2
+7\cdot\jj 212222^2-4\cdot\jj 222222\cdot\jj 212222+1
\end{align*}
and
\begin{align*}
  \mathrm{Nf}_>&\left(\widetilde{tv}_{2,1}\left(L(8,3)\#L(7,2)\right)\right)\\ &=w(2)^3\cdot\jj 222222^2+\jj 212212^2 +
  9\cdot\jj 212222^2
  -6\cdot\jj 222222\cdot\jj 212222 +1\\
  \not=\mathrm{Nf}_>&\left(\mathrm{Nf}_>\left(\widetilde{tv}_{2,1}\left(L(8,3)\right)\right)
    \cdot \mathrm{Nf}_>\left(\widetilde{tv}_{2,1}\left(L(7,2)\right)\right)\right) \\
  &= w(2)^4\cdot\jj 222222^2+2\cdot w(2)^3\cdot\jj 222222^2
-\jj 212222^2+\jj 222222\cdot\jj 212222+1
\end{align*}
For our computations, we used the following special spines (in Matveev-coding):
\begin{itemize}
\item (( 1, 1, 2, -3), ( 1, 3, -4, -4, -2), ( 2, -4, -3)) for
  $L(7,2)$,
\item ((1, 1, 2, -3), ( 1, 3, -4, -2), ( 2, -4, -4, -3)) for $L(8,3)$, 
\item ((1, 1, 2, -3), (1, 3, -7, -16, -15, -4, -2), (2, -7, -6, -5,
  -4, -7, 17, 18, -4, -3), (8, 8, 9, -10), (8, 10, -14, -17, -16, -11,
  -14, 18, 15, -11, -9), (9, -14, -13, -12, -11, -10), (15, 12, -5),
  (16, -6, -12), (17, -13, 6), (18, 5, 13)) for $L(8,3)\#L(7,2)$, and
\item ((1, 1, 2, -3), (1, 3, -7, -16, -15, -4, -2), (2, -7, -6, -5,
  -4, -7, 17, 18, -4, -3), (8, 8, 9, -10), (8, 10, -14, -13, -12, -11,
  -9), (9, -14, -17, -16, -11, -14, 18, 15, -11, -10), (15, 12, -5),
  (16, -6, -12), (17, -13, 6), (18, 5, 13)) for $L(8,3)\#L(8,3)$.
\end{itemize}
\end{proof}

The first statement of Proposition~\ref{thm:results} says, in
combination with Theorem~\ref{thm:NumVsId}, that one should expect
that ideal Turaev--Viro invariants are, in general, stronger than a
combination of all associated numerical Turaev--Viro invariants.

The second and third statement of Proposition~\ref{thm:results} shows
that $\widetilde {tv}_{2,1}(\cdot)$ sees properties of manifolds that
are invisible for homology, and vice versa.

Statement~(4) is surprising, because one would expect that one obtains
a stronger invariant if one avoids to impose simplifying assumptions.
But this is not necessarily the case. Statement~(5) is even more
surprising, because by statement~(1) the Turaev--Viro ideals are not
radical --- hence there are elements of $R$ so that the cosets with
respect to $\sqrt{I_{m,n}}$ coincide, but the cosets with respect to
$I_{m,n}$ are different. Are there compact $3$-manifolds $M_1$, $M_2$
that can be distinguished by some ideal Turaev--Viro invariant
$tv(\cdot)$ but can not be distinguished by all associated numerical
Turaev--Viro invariants, i.e., can not be distinguished by
$\widehat{tv}(\cdot)$? Note that $\widetilde{tv}_{2,1}$ is stronger
than the $\eps$-invariant; but the $\eps$-invariant is not the only
numerical Turaev-Viro invariant associated to $\widetilde{tv}_{2,1}$
(see~\cite[Sec.~8.1]{matveev}).

The last statement of Proposition~\ref{thm:results} is a bad news if
one wants to construct a Topological Quantum Field Theory. But it is a
good news if one aims to construct invariants that potentially detect
counter-examples of the Andrews-Curtis conjecture. Namely, by a result
of Bobtcheva and Quinn~\cite{BobtchevaQuinn}, an invariant for
Andrews-Curtis moves descending from a multiplicative invariant of
$4$-thickenings of special $2$-polyhedra only depends on homology if
the Euler characteristic of the $2$-complex under consideration is at
least $1$. But a non-multiplicative ideal Turaev--Viro invariant for
Andrews-Curtis moves~\cite{kingAC} is potentially more useful. Note
that the $\eps$-invariant is multiplicative, and is $1+\eps$ on both
$L(8,3)$ and $L(7,2)$.

\subsection*{Acknowledgement}
I'm grateful to S. V. Matveev for providing me with a list of one
special spine for each closed orientable irreducible $3$-manifold up
to complexity $9$, which was essential for the computations. I thank
I. Bobtcheva for interesting discussions.

\vspace*{1cm}
\noindent
\textit{Simon A. King\\
Department of Mathematics\\
Technical University Darmstadt\\
64289 Darmstadt\\
GERMANY}\\
\texttt{king@mathematik.tu-darmstadt.de}

\end{document}